\documentclass[12pt, reqno]{amsart}
\usepackage{amssymb, amsmath, amsfonts, amscd, calligra, mathrsfs}
\usepackage[raiselinks=false,colorlinks=true,citecolor=blue,urlcolor=blue,
linkcolor=blue,bookmarksopen=true,pdftex]{hyperref}
\usepackage[dvipsnames]{xcolor}
\usepackage{tikz}\usepackage{tikz-cd}
\newtheorem{theorem}{Theorem}

\newtheorem{lemma}[theorem]{Lemma}
\newtheorem{remark}[theorem]{Remark}

\newcommand{\rank}{\textrm{rank}}
\begin{document}
\baselineskip=14.5pt
\title{HOMOLOGICAL METHODS IN CERTAIN PICARD GROUP COMPUTATIONS}
\author[P. Biswas]{Pritthijit Biswas}
\address{Chennai Mathematical Institute, H1 Sipcot IT Park, Siruseri, Kelambakkam 603103.}
\email{pritthijit@cmi.ac.in}
\begin{abstract}
Let $G$ be a connected complex semisimple Lie group, $\Gamma$ be a cocompact, irreducible and torsionless lattice in $G$ and $K$ be a maximal compact subgroup of $G$. Assume $\Gamma$ acts by left multiplication and $K$ acts by right multiplication on $G$. Let $M_{\Gamma}= \Gamma\backslash G$, $X=G/K$ and $X_{\Gamma}=\Gamma\backslash X$. In this article we prove that for any $n\geq0$, the composition $H^{n}(X_{\Gamma},\mathbb{C})\rightarrow H^{n}(M_{\Gamma},\mathbb{C})\rightarrow H^{n}(M_{\Gamma},\mathcal{O}_{M_{\Gamma}})$ is an isomorphism. As an application when $G$ is simply connected, we compute the Picard group of $M_{\Gamma}$ for the cases rank($G$) $=1,2$. More precisely we show that if rank($G$) $=1$, $Pic(M_{\Gamma})=(\mathbb{C}^{r}/\mathbb{Z}^{r})\oplus A$ and if rank($G$) $=2$, then $Pic(M_{\Gamma})\cong A$ via the first Chern class map, where $A$ is the torsion subgroup of $H^{2}(M_{\Gamma},\mathbb{Z})$ and $r$ is the rank of $\Gamma/[\Gamma,\Gamma]$. 
\end{abstract} 
\maketitle
\section{Introduction}
Recall that for any connected complex manifold $M$, its Picard group denoted by $Pic(M)$ is the set of isomorphism classes of holomorphic line bundles on $M$ which forms an abelian group with respect to the tensor product operation. $Pic(M)$ is naturally identified with $H^{1}(M,\mathcal{O}^{*}_{M})$ and $Pic^{0}(M)$ is defined to be the kernel of the first Chern class map $c_{1}:Pic(M)\rightarrow H^{2}(M,\mathbb{Z})$. See \cite{huybrechts} for further details. 

Let $G$ be a semisimple Lie group with finitely many connected components. A discrete subgroup $\Gamma$ of $G$ is called a {\it lattice} in $G$ if $\Gamma\backslash G$ ($\Gamma$ assumed to act on $G$ by left multiplication)  has a right $G$-invariant measure (i.e. an invariant measure for the right action of $G$ on $\Gamma\backslash G$) with finite total mass. If $\Gamma\backslash G$ is compact, then $\Gamma$ is said to be {\it cocompact} (or {\it uniform}). It is well known due to Borel and Harish-Chandra, that any real semisimple Lie group with finitely many connected components has both cocompact and non-cocompact lattices. A lattice $\Gamma$ in $G$ is said to be {\it reducible} if there exists a finite sheeted covering group $\pi:\tilde G\to G$ such that $\tilde G=\tilde G_1\times \tilde G_2$ where $\tilde G_j$ are non-compact subgroups of $\tilde G$ and a subgroup $\Lambda\subset \Gamma$ having finite index in $\Gamma$ such that $\pi^{-1}(\Lambda) = \Lambda_1\times \Lambda_2$ where $\Lambda_j\subset \tilde G_j$ are lattices. If $\Gamma$ is not reducible, then it is said to be {\it irreducible}. We refer the reader to the book by Raghunathan \cite{raghunathan} for the definition of lattices in any general Lie group, basic theorems about them and equivalent formulations of irreducible lattices.   

Recall that the real rank of a linear connected real semisimple Lie group $G$, denoted $\rank_\mathbb R(G)$ is the maximum dimension  $r$ of a diagonalizable subgroup of $G$ isomorphic to $(\mathbb R_{>0})^r $. If $G$ is a connected complex semisimple Lie group, we define its real rank by regarding it as a real Lie group by forgetting the 
complex structure.  In this case, the real rank of $G$ equals the rank of $G$ viewed as a complex algebraic group, namely, the dimension of a maximal (algebraic) torus (isomorphic to $(\mathbb C^*)^r$) contained in $G$. For example, when $G=SL(n,\mathbb C)$, we have $\rank_\mathbb RG=n-1$.  In view of this, we shall henceforth write $\rank(G)$ to mean $\rank_\mathbb R(G)$.  

Let $G$ be a connected complex semisimple Lie group, $\Gamma$ be a cocompact, irreducible and torsionless lattice in $G$ and $K$ be a maximal compact subgroup of $G$. Assume $\Gamma$ acts by left multiplication and $K$ acts by right multiplication on $G$. Let $M_{\Gamma}= \Gamma\backslash G$, $X=G/K$ and $X_{\Gamma}=\Gamma\backslash X$. Note that $X$, $X_{\Gamma}$ are smooth manifolds and does not necessarily admit complex structures. It is well known that $X$ is contractible. We refer the reader to the texts \cite{V}, \cite{L} and \cite{K} for basic facts on complex semisimple Lie groups including their structure theory and classification. Since $X\rightarrow X_{\Gamma}$ is a covering projection, $X_{\Gamma}=K(\Gamma,1)$. We obtain the following result.
\begin{theorem}\label{t1}
Let $G$ be a connected complex semisimple Lie group, $\Gamma$ be a cocompact, irreducible and torsionless lattice in $G$ and $K$ be a maximal compact subgroup of $G$. Assume $\Gamma$ acts by left multiplication and $K$ acts by right multiplication on $G$. Then for any $n\geq 0$, the composition $H^{n}(X_{\Gamma},\mathbb{C})\rightarrow H^{n}(M_{\Gamma},\mathbb{C})\rightarrow H^{n}(M_{\Gamma},\mathcal{O}_{M_{\Gamma}})$ is an isomorphism.
\end{theorem} 

For $p,q\geq 0$, see \cite{huybrechts} for the definition of Dolbeault cohomologies $H^{p,q}(M)$ for any complex manifold $M$. We prove theorem \ref{t1} by showing that for all $n\geq 0$, the composition $H^{n}(\Gamma, \mathbb{C})\cong H^{n}(X_{\Gamma},\mathbb{C})\rightarrow H^{n}(M_{\Gamma},\mathbb{C})\rightarrow H^{n}(M_{\Gamma},\mathcal{O}_{M_{\Gamma}})$ is same as the isomorphism $H^{n}(\Gamma, \mathbb{C})\rightarrow H^{n}(M_{\Gamma},\mathcal{O}_{M_{\Gamma}})$ obtained by Akhiezer in \cite[Theorem 1]{akhiezer}. However the proof of theorem \ref{t1} is approached after formulating the problem in a more general setting and involves multiple use of the spectral sequences due to Grothendieck as given in \cite[Theorem 2.4.1]{grothendieck} and finding key naturality correspondences between certain initial terms as would be formally elucidated in the sequel. As an application of theorem \ref{t1}, we obtain the following result, which has been the motivation to prove theorem \ref{t1}.
\begin{theorem}\label{t2}
Let $G$ be a simply connected complex semisimple Lie group, $\Gamma$ be a cocompact, irreducible and torsionless lattice in $G$ which acts by left multiplication on $G$. Then  if rank($G$) $=1$, $Pic(M_{\Gamma})=(\mathbb{C}^{r}/\mathbb{Z}^{r})\oplus A$ and if rank($G$) $=2$, $Pic(M_{\Gamma})\cong A$ via the first Chern class map, where $A$ is the torsion subgroup of $H^{2}(M_{\Gamma},\mathbb{Z})$ and $r$ is the rank of $\Gamma/[\Gamma,\Gamma]$.
\end{theorem}
Considering the hypotheses of theorem \ref{t2}, $Pic^{0}(M_{\Gamma})$ for any rank of $G$ and $Pic(M_{\Gamma})$ for rank($G$) $\geq 3$ have been computed by Biswas and Sankaran in \cite[Theorem 1]{picardgroups}. 
\section{General Formulation}
For the rest of this article, we shall use the notion of proper Lie group action on a smooth manifold as given in \cite[\S 21]{lee}. Let $Y$ be a connected complex manifold and $\Gamma$ be a discrete group which acts holomorphically freely and properly on the left of $Y$ (hence is a covering space action). Let $M_{\Gamma}$ denote the complex manifold $\Gamma\backslash Y$. Let $K$ be a connected real Lie group which acts smoothly, freely and properly on the right of $Y$ such that $\gamma.(y.k)=(\gamma.y).k$ $\forall$ $\gamma\in \Gamma$, $y\in Y$ and $k\in K$. Let $X$ denote the smooth manifold $Y/K$. Then it is clear that $\Gamma$ acts continuously on the left of $X$. Moreover assume that $H^{q}(X,\mathbb{C})=0$ for all $q\gneq 0$. Let $X_{\Gamma}$ denote the topological space $\Gamma\backslash X$.

Let $\pi:Y\rightarrow M_{\Gamma}$ be the usual projection which is a holomorphic covering. Let $p:Y\rightarrow X$, $\tilde{\pi}:X\rightarrow X_{\Gamma}$ be the usual projections ($p$ is smooth). There exists a unique open surjection (also a quotient map) $\tilde{p}:M_{\Gamma}\rightarrow X_{\Gamma}$, such that $\tilde{p}\circ \pi=\tilde{\pi}\circ p$. Henceforth we consider the following notions :
\begin{itemize}
\item For any topological space $Z$ and a continuous left action of a topological group $G$ on $Z$, let $C^{Z(G)}$ denote the abelian category of $G$-sheaves of abelian groups on $Z$. It is well known that $C^{Z(G)}$ has enough injectives \cite[Proposition 5.1.2]{grothendieck}. For any topological space $Z$ we prefer the notation $C^{Z}$ for the abelian category of sheaves of abelian groups on $Z$. Let $\Lambda_{Z} : C^{Z(G)}\rightarrow G$-Mod denote the usual global section functor, where $G$-Mod denotes the abelian category of left $G$-modules. Let $\Lambda^{G} : C^{Z(G)}\rightarrow C^{Z}$ denote the $G$-fixed point sheaf functor, when $G$ acts trivially on $Z$. Let $\Lambda_{Z}:C^{Z}\rightarrow Ab$ denote the usual global section functor, where $Ab$ denotes the abelian category of abelian groups.
\item We consider the spaces $M_{\Gamma}$, $X_{\Gamma}$ to have the trivial left action of $\Gamma$. 
\item For any topological space $Z$, let $\underline{\mathbb{C}}_{Z}$ denote the constant sheaf $\mathbb{C}$. We will treat $\underline{\mathbb{C}}_{X}$ and $\underline{\mathbb{C}}_{Y}$ as objects of $C^{X(\Gamma)}$ and $C^{Y(\Gamma)}$ (via composition with left multiplications by elements of $\Gamma$) respectively, and will treat $\underline{\mathbb{C}}_{X_{\Gamma}}$ and $\underline{\mathbb{C}}_{M_{\Gamma}}$ as objects of $C^{X_{\Gamma}}$ and $C^{M_{\Gamma}}$ respectively. For any complex manifold $Z$, let $\mathcal{O}_{Z}$ denote its structure sheaf. We use the notation $\mathscr{M}$ for the abelian category of left $\Gamma$-modules.
\end{itemize}
For any $n\geq 0$, the goal is to identify the composition ${\iota_{*}\circ \tilde{p}^{*}} : H^{n}(X_{\Gamma},\mathbb{C})\overset{\tilde{p}^{*}}\rightarrow H^{n}(M_{\Gamma},\mathbb{C})\overset{\iota_{*}}\rightarrow H^{n}(M_{\Gamma},\mathcal{O}_{M_{\Gamma}})$ where $\iota$ is the inclusion $\underline{\mathbb{C}}_{M_{\Gamma}}\hookrightarrow \mathcal{O}_{M_{\Gamma}}$, with a homomorphism $\psi$, where $\psi$ is obtained using a suitable network of functors between certain abelian categories of sheaves. This involves multiple use of the Grothendieck spectral sequences (\cite[Theorem 2.4.1]{grothendieck}). In particular when we are in the setting of theorem \ref{t1}, the map $\psi$ coincides with the Akhiezer's isomorphism $H^{n}(\Gamma, \mathbb{C})\rightarrow H^{n}(M_{\Gamma},\mathcal{O}_{M_{\Gamma}})$ as given in \cite[Theorem 1]{akhiezer} when $p=0$. This needs some work and will be made rigorous as we go along. We start with :

\usetikzlibrary{matrix,calc}
\begin{equation}
\begin{tikzpicture}[-stealth,
  label/.style = { font=\footnotesize }, baseline=(current  bounding  box.center)]
  \matrix (m)
    [
      matrix of math nodes,
      row sep    = 4.5em,
      column sep = 3.5em
    ]
    {
      \mathscr{M} & C^{Y(\Gamma)} & C^{M_{\Gamma}(\Gamma)} & C^{M_{\Gamma}} & Ab\\
      \mathscr{M} & C^{X(\Gamma)} & C^{X_{\Gamma}(\Gamma)} & C^{X_{\Gamma}} & Ab \\
    };
  \path (m-1-2) edge node [above, label] {$\pi_{*}$} (m-1-3);
  \path (m-2-2) edge node [above, label] {$\tilde{\pi}_{*}$} (m-2-3);
  \path (m-1-5) edge node [left, label] {$id$} (m-2-5);
  \path (m-1-5) edge node [right, label] {$\cong$} (m-2-5);
  \path (m-1-1) edge node [left, label] {$\cong$} (m-2-1);
  \path (m-1-1) edge node [right, label] {$id$} (m-2-1);
  \path (m-1-2) edge node [right, label] {$p_{*}$} (m-2-2);
  \path (m-1-3) edge node [above, label] {$\Lambda^{\Gamma}$} (m-1-4);
  \path (m-2-3) edge node [above, label] {$\Lambda^{\Gamma}$} (m-2-4);
  \path (m-1-4) edge node [right, label] {$\tilde{p}_{*}$} (m-2-4);
  \path (m-1-4) edge node [above, label] {$\Lambda_{M_{\Gamma}}$} (m-1-5);
  \path (m-2-4) edge node [above, label] {$\Lambda_{X_{\Gamma}}$} (m-2-5);
  \path (m-1-2) edge [out=30, in=150] node [above, label] {$\pi_{*}^{\Gamma}$} (m-1-4);
  \path (m-1-1) edge [out=45, in=135] node [above, label] {$F^{\Gamma}$} (m-1-5);
  \path (m-2-1) edge [out=-45, in=-135] node [below, label] {$F^{\Gamma}$} (m-2-5);
  \path (m-2-2) edge [out=-30, in=-150] node [below, label] {$\tilde{\pi}_{*}^{\Gamma}$} (m-2-4);
  \path (m-1-2) edge node [above, label] {$\Lambda_{Y}$} (m-1-1);
  \path (m-2-2) edge node [above, label] {$\Lambda_{X}$} (m-2-1);
\end{tikzpicture}
\end{equation}
In (1), $p_{*}$, $\tilde{p}_{*}$, $\pi_{*}$, $\tilde{\pi}_{*}$ denote the usual direct image sheaf functors and $F^{\Gamma}$ denotes the usual $\Gamma$-fixed point module functor. $\pi^{\Gamma}_{*}:= \Lambda^{\Gamma}\circ \pi_{*}$ and $\tilde{\pi}^{\Gamma}_{*}:= \Lambda^{\Gamma}\circ \tilde{\pi}_{*}$ and all the functors in (1) are left exact, covariant and additive. One can directly verify the commutativity of (1) upto equality of functors. 
Straightforward arguments yield that $p_{*}$, $\tilde{p}_{*}$, $\pi_{*}$, $\tilde{\pi}_{*}$ take injectives to injectives. 

$\pi^{\Gamma}_{*}$, $\tilde{\pi}^{\Gamma}_{*}$ take injectives to flasques and hence acyclics (with respect to the global section functors respectively in each case). $\Lambda_{Y} : C^{Y(\Gamma)}\rightarrow \mathscr{M}$ and $\Lambda_{X} : C^{X(\Gamma)}\rightarrow \mathscr{M}$ take injectives to injectives (\cite[Ch. 5, Corollary after Proposition 5.1.3]{grothendieck}).

Note that there exists an isomorphism $\underline{\mathbb{C}}_{X}{\overset{t}{\underset{\cong}\rightarrow}} p_{*}\underline{\mathbb{C}}_{Y}$ in $C^{X(\Gamma)}$ and hence an isomorphism $\underline{\mathbb{C}}_{X}(X)\overset{\Lambda_{X}(t)}{\underset{\cong}\rightarrow} p_{*}\underline{\mathbb{C}}_{Y}(X)=\underline{\mathbb{C}}_{Y}(Y)$ in $\mathscr{M}$. It is easy to see that $\pi^{\Gamma}_{*}\underline{\mathbb{C}}_{Y}\cong\underline {\mathbb{C}}_{M_{\Gamma}}$, $\pi^{\Gamma}_{*}\mathcal{O}_{Y}\cong \mathcal{O}_{M_{\Gamma}}$ in $C^{M_{\Gamma}}$ and $\tilde{\pi}^{\Gamma}_{*}\underline{\mathbb{C}}_{X}\cong \underline{\mathbb{C}}_{X_{\Gamma}}$ in $C^{X_{\Gamma}}$. From the commutativity of (1) we get that  $\tilde{\pi}^{\Gamma}_{*}p_{*}{\underline{\mathbb{C}}}_{Y}\cong\tilde{p}_{*}{\underline{\mathbb{C}}}_{M_{\Gamma}}$. There is an obvious isomorphism $H^{n}(\Gamma, \underline{\mathbb{C}}_{X}(X))\overset{\Lambda_{X}(t)_{*}}{\underset{\cong}\longrightarrow}H^{n}(\Gamma, p_{*}\underline{\mathbb{C}}_{Y}(X))$ and we have the natural map $\Lambda_{Y}(\iota)_{*}: H^{n}(\Gamma, {\underline{\mathbb{C}}}_{Y}(Y))\rightarrow H^{n}(\Gamma, \mathcal{O}_{Y}(Y))$, where $\iota: \underline{\mathbb{C}}_{Y}\hookrightarrow \mathcal{O}_{Y}$ is the usual inclusion of $\Gamma$-sheaves of abelian groups on $Y$. Let $n\ge 0$ be arbitrary.
The Grothendieck Spectral Sequence (GSS) associated to the setup $C^{X(\Gamma)}\overset{\tilde{\pi}^{\Gamma}_{*}}\rightarrow C^{X_{\Gamma}}\overset{\Lambda_{X_{\Gamma}}}\rightarrow Ab$  gives the edge homomorphism (by which, for the rest of the article, we will always mean $E^{n,0}_{2}\rightarrow E^{n,0}_{\infty}$ with respect to the row wise filtration of the associated double complex used to define it unless otherwise stated, see \cite[Theorem 10.47]{rotman}) $H^{n}(X_{\Gamma}, \tilde{\pi}^{\Gamma}_{*}\underline{\mathbb{C}}_{X})\overset{s}\rightarrow R^{n}(\Lambda_{X_{\Gamma}}\tilde{\pi}^{\Gamma}_{*})(\underline{\mathbb{C}}_{X})$. GSS associated to the setup $C^{X(\Gamma)}\overset{\Lambda_{X}}\rightarrow \mathscr{M}\overset{F^{\Gamma}}\rightarrow Ab$ gives the edge homomorphism $H^{n}(\Gamma,\underline{\mathbb{C}}_{X}(X))\overset{u}{\underset{\cong}\rightarrow} R^{n}(\Lambda_{X_{\Gamma}}\tilde{\pi}^{\Gamma}_{*})(\underline{\mathbb{C}}_{X})$. By our hypothesis we have $H^{q}(X,\mathbb{C})=0$ for all $q>0$, which implies that $u$ is an isomorphism \cite [Ch. 5, Proposition 5.2.5]{grothendieck}. GSS associated to the setup $C^{Y(\Gamma)}\overset{\Lambda_{Y}}\rightarrow \mathscr{M}\overset{F^{\Gamma}}\rightarrow Ab$ gives the edge homomorphism $H^{n}(\Gamma, \mathcal{O}_{Y}(Y))\overset{v}\rightarrow R^{n}(\Lambda_{M_{\Gamma}}\pi^{\Gamma}_{*})(\mathcal{O}_{Y})$. GSS associated to the setup $C^{Y(\Gamma)}\overset{\pi^{\Gamma}_{*}}\rightarrow C^{M_{\Gamma}}\overset{\Lambda_{M_{\Gamma}}}\rightarrow Ab$ gives the edge homomorphism $H^{n}(M_{\Gamma}, \mathcal{O}_{M_{\Gamma}})\overset{w}{\underset{\cong}\rightarrow} R^{n}(\Lambda_{M_{\Gamma}}\pi^{\Gamma}_{*})(\mathcal{O}_{Y})$. Since $\pi$ is a holomorphic covering, $w$ is an isomorphism \cite [Ch. 5, Corollary 1 after Theorem 5.3.1]{grothendieck}.

Keeping all the above notions in mind, we can now finally state the central theorem of this article :
\begin{theorem}\label{theorem1} After identifying $H^{n}(X_{\Gamma}, \tilde{\pi}^{\Gamma}_{*}\underline{\mathbb{C}}_{X})$ with the singular cohomology of $X_{\Gamma}$ with complex coefficients - $H^{n}(X_{\Gamma},\mathbb{C})$,  the map $\psi$ defined by the composite $w^{-1}\circ v\circ \Lambda_{Y}(\iota)_{*}\circ \Lambda_{X}(t)_{*}\circ u^{-1} \circ s$ makes the following diagram commutative :
\usetikzlibrary{matrix,calc}
\begin{equation}
\begin{tikzpicture}[-stealth,
  label/.style = { font=\footnotesize }, baseline=(current  bounding  box.center)]
  \matrix (m)
    [
      matrix of math nodes,
      row sep    = 4.5em,
      column sep = 4em
    ]
    {
      H^{n}(X_{\Gamma}, \mathbb{C}) & H^{n}(M_{\Gamma}, \mathcal{O}_{M_{\Gamma}})\\
      H^{n}(M_{\Gamma}, \mathbb{C}) & H^{n}(M_{\Gamma}, \mathcal{O}_{M_{\Gamma}}) \\
    };
    \path (m-1-1) edge node [above, label] {$\psi$} (m-1-2);
    \path (m-1-1) edge node [right, label] {$\tilde{p}^{*}$} (m-2-1);
    \path (m-2-1) edge node [above, label] {$\iota_{*}$} (m-2-2);
    \path (m-1-2) edge node [right, label] {$id$} (m-2-2);
    \path (m-1-2) edge node [left, label] {$\cong$} (m-2-2);
\end{tikzpicture}
\end{equation}   
\end{theorem}
\section{Proofs of Theorems 1 and 3}
In this section we first provide a proof of theorem \ref{theorem1} and we use it to prove theorem \ref{t1}. Unless otherwise mentioned, for the rest of this section, we assume that $Y$ is a connected complex manifold, $\Gamma$ is a discrete group acting holomorphically, freely and properly on the left of $Y$. Let $K$ be a connected real Lie group which acts smoothly, freely and properly on the right of $Y$ such that $\gamma.(y.k)=(\gamma.y).k$ $\forall$ $\gamma\in \Gamma$, $y\in Y$ and $k\in K$. As before define $M_{\Gamma}:=\Gamma\backslash Y$, $X:= Y/K$. Clearly $\Gamma$ acts on the left of $X$. Let $X_{\Gamma}$ be the topological space $\Gamma\backslash X$. As before we have commutativity of (1) and the natural isomorphisms $\pi^{\Gamma}_{*}\underline{\mathbb{C}}_{Y}\cong\underline {\mathbb{C}}_{M_{\Gamma}}$, $\pi^{\Gamma}_{*}\mathcal{O}_{Y}\cong \mathcal{O}_{M_{\Gamma}}$ and $\tilde{\pi}^{\Gamma}_{*}\underline{\mathbb{C}}_{X}\cong\underline {\mathbb{C}}_{X_{\Gamma}}$. Let $n\geq 0$. 
GSS associated to the setup $C^{X(\Gamma)}\overset{\tilde{\pi}^{\Gamma}_{*}}\rightarrow C^{X_{\Gamma}}\overset{\Lambda_{X_{\Gamma}}}\rightarrow Ab$, gives the edge homomorphism $H^{n}(X_{\Gamma}, \tilde{\pi}^{\Gamma}_{*}p_{*}\underline{\mathbb{C}}_{Y}\cong\tilde{p}_{*}\underline{\mathbb{C}}_{M_{\Gamma}})\overset{e_{1}}\rightarrow R^{n}(\Lambda_{X_{\Gamma}}\tilde{\pi}^{\Gamma}_{*})(p_{*}\underline{\mathbb{C}}_{Y})$. Similarly GSS associated to the setup $C^{Y(\Gamma)}\overset{\pi^{\Gamma}_{*}}\rightarrow C^{M_{\Gamma}}\overset{\Lambda_{M_{\Gamma}}}\rightarrow Ab$, gives the edge homomorphism $H^{n}(M_{\Gamma}, \pi^{\Gamma}_{*}\underline{\mathbb{C}}_{Y})\overset{e_{2}}\rightarrow R^{n}(\Lambda_{M_{\Gamma}}\pi^{\Gamma}_{*})(\underline{\mathbb{C}}_{Y})$.
\begin{lemma}\label{lemma2}
We keep the above notations. There exists a map $\mu : R^{n}(\Lambda_{X_{\Gamma}}\tilde{\pi}^{\Gamma}_{*})(p_{*}\underline{\mathbb{C}}_{Y})\rightarrow R^{n}(\Lambda_{M_{\Gamma}}\pi^{\Gamma}_{*})(\underline{\mathbb{C}}_{Y})$ which makes the following diagram commutative :
\begin{equation}
\begin{tikzpicture}[-stealth,
  label/.style = { font=\footnotesize }, baseline=(current  bounding  box.center)]
  \matrix (m)
    [
      matrix of math nodes,
      row sep    = 4.5em,
      column sep = 3em
    ]
    {
      R^{n}(\Lambda_{X_{\Gamma}}\tilde{\pi}^{\Gamma}_{*})(p_{*}\underline{\mathbb{C}}_{Y}) & R^{n}(\Lambda_{M_{\Gamma}}\pi^{\Gamma}_{*})(\underline{\mathbb{C}}_{Y})\\
      H^{n}(X_{\Gamma},\tilde{p}_{*}\underline{\mathbb{C}}_{M_{\Gamma}}) & H^{n}(M_{\Gamma},\underline{\mathbb{C}}_{M_{\Gamma}}) \\
    };
    \path (m-1-1) edge node [above, label] {$\mu$} (m-1-2);
    \path (m-2-1) edge node [right, label] {$e_{1}$} (m-1-1);
    \path (m-2-1) edge node [above, label] {via $\tilde{p}$} (m-2-2);
    \path (m-2-2) edge node [right, label] {$e_{2}$} (m-1-2);
\end{tikzpicture}
\end{equation}   
\end{lemma}
{\it Proof.}
Choose injective resolutions $0\rightarrow {\underline{\mathbb{C}}}_{Y}\rightarrow \mathscr{I}^{\bullet}$ in $C^{Y(\Gamma)}$ and $0\rightarrow p_{*}{{\underline{\mathbb{C}}}_{Y}}\rightarrow \mathscr{J}^{\bullet}$ in $C^{X(\Gamma)}$. Then $0\rightarrow \tilde{\pi}^{\Gamma}_{*}p_{*}{\underline{\mathbb{C}}}_{Y}\rightarrow \tilde{\pi}^{\Gamma}_{*}\mathscr{J}^{\bullet}$ is a complex in $C^{X_{\Gamma}}$. Let $0\rightarrow \tilde{\pi}^{\Gamma}_{*}\mathscr{J}^{\bullet}\rightarrow \mathscr{J}^{\bullet, *}$ be a Cartan-Eilenberg resolution (see \cite[Ch. 10]{rotman} for definition) in $C^{X_{\Gamma}}$. Let $\tilde{p}^{-1}:C^{X_{\Gamma}}\rightarrow C^{M_{\Gamma}}$ and $p^{-1}:C^{X(\Gamma)}\rightarrow C^{Y(\Gamma)}$ be the usual inverse image sheaf functors.

There is an obvious natural transformation $\zeta:{\tilde{p}}^{-1}{\tilde{\pi}}^{\Gamma}_{*}\rightarrow \pi^{\Gamma}_{*}p^{-1}$, which gives us a map of complexes $\zeta_{\mathscr{J}^{\bullet}}:{\tilde{p}}^{-1}{\tilde{\pi}}^{\Gamma}_{*}\mathscr{J}^{\bullet}\rightarrow \pi^{\Gamma}_{*}(p^{-1}\mathscr{J}^{\bullet})$ extending the morphism $\zeta_{p_{*}{\underline{\mathbb{C}}}_{Y}}:{\tilde{p}}^{-1}{\tilde{\pi}}^{\Gamma}_{*}(p_{*}{\underline{\mathbb{C}}}_{Y})\rightarrow \pi^{\Gamma}_{*}(p^{-1}p_{*}{\underline{\mathbb{C}}}_{Y})$. Since $p^{-1}$ is exact, there exists a morphism of complexes $p^{-1}\mathscr{J}^{\bullet}\overset{c^{\bullet}}\rightarrow \mathscr{I}^{\bullet}$ extending the obvious morphism $p^{-1}p_{*}{\underline{\mathbb{C}}}_{Y}\overset{c}\rightarrow \underline{\mathbb{C}}_{Y}$. 
So we have a map of complexes $\pi^{\Gamma}_{*}(p^{-1}\mathscr{J}^{\bullet})\overset{\pi^{\Gamma}_{*}(c^{\bullet})}\longrightarrow \pi^{\Gamma}_{*}\mathscr{I}^{\bullet}$ extending the morphism $\pi^{\Gamma}_{*}(p^{-1}p_{*}{\underline{\mathbb{C}}}_{Y})\overset{\pi^{\Gamma}_{*}(c)}\longrightarrow \pi^{\Gamma}_{*}{\underline{\mathbb{C}}}_{Y}$. Thus we get a morphism of complexes ${\tilde{p}}^{-1}{\tilde{\pi}}^{\Gamma}_{*}\mathscr{J}^{\bullet}\overset{\delta^{\bullet}}\rightarrow \pi^{\Gamma}_{*}\mathscr{I}^{\bullet}$ which extends the morphism ${\tilde{p}}^{-1}{\tilde{\pi}}^{\Gamma}_{*}(p_{*}{\underline{\mathbb{C}}}_{Y})\overset{\delta}\rightarrow \pi^{\Gamma}_{*}{\underline{\mathbb{C}}}_{Y}$, where $\delta$ := $\pi^{\Gamma}_{*}(c)\circ \zeta_{p_{*}\underline{\mathbb{C}}_{Y}}$ and $\delta^{\bullet}$ := $\pi^{\Gamma}_{*}(c^{\bullet})\circ \zeta_{\mathscr{J}^{\bullet}}$. We choose a Cartan-Eilenberg resolution $0\rightarrow \pi^{\Gamma}_{*}\mathscr{I}^{\bullet}\rightarrow \mathscr{I}^{\bullet, *}$ of the complex $0\rightarrow \pi^{\Gamma}_{*}{\underline{\mathbb{C}}}_{Y}\rightarrow \pi^{\Gamma}_{*}\mathscr{I}^{\bullet}$ in $C^{M_{\Gamma}}$. Note that we also have a double complex with exact columns namely $0\rightarrow \tilde{p}^{-1}{\tilde{\pi}}^{\Gamma}_{*}\mathscr{J}^{\bullet}\rightarrow \tilde{p}^{-1}\mathscr{J}^{\bullet, *}$. 
Using the explicit construction of Cartan-Eilenberg resolution of a complex (see \cite[Theorem 10.45]{rotman} and its proof), standard homological algebra techniques and the fact that ${\tilde{p}}^{-1}$ is an exact functor which commutes with direct sums, we get a map of double complexes $\tilde{p}^{-1}\mathscr{J}^{\bullet, *}\rightarrow \mathscr{I}^{\bullet, *}$ extending the map of complexes $\delta^{\bullet} : {\tilde{p}}^{-1}{\tilde{\pi}}^{\Gamma}_{*}\mathscr{J}^{\bullet}\rightarrow \pi^{\Gamma}_{*}\mathscr{I}^{\bullet}$ in $C^{M_{\Gamma}}$. We also have a natural transformation $\eta: \Lambda_{X_{\Gamma}}\rightarrow \Lambda_{M_{\Gamma}}{{\tilde{p}}^{-1}}$. Thus we get the following commutative diagram :

\usetikzlibrary{matrix,calc}
\vspace{-1em}
{\tiny
\begin{equation}
\begin{tikzpicture}[-stealth,
  label/.style = { font=\footnotesize }, baseline=(current  bounding  box.center)]
  \matrix (m)
    [
      matrix of math nodes,
      row sep    = 4em,
      column sep = 4.5em
    ]
   {
{^{(r)}}{E^{n,0}_{2}}=H^{n}(M_{\Gamma}, \pi^{\Gamma}_{*}{{\underline{\mathbb{C}}}_{Y}}) & H^{n}(tot(\Lambda_{M_{\Gamma}}{\mathscr{I}}^{\bullet, *})) & ^{(c)}E^{n,0}_{2}=R^{n}(\Lambda_{M_{\Gamma}}\pi^{\Gamma}_{*})({\underline{\mathbb{C}}}_{Y}) \\
{^{(r)}}{{\hat{E}}^{n,0}_{2}} & H^{n}(tot(\Lambda_{M_{\Gamma}}{\tilde{p}}^{-1}\mathscr{J}^{\bullet, *})) & ^{(c)}{{\hat{E}}^{n,0}_{2}}=H^{n}(\Lambda_{M_{\Gamma}}{\tilde{p}}^{-1}{\tilde{\pi}}^{\Gamma}_{*}\mathscr{J}^{\bullet})\\
{^{(r)}}{{\tilde{E}}^{n,0}_{2}}=H^{n}(X_{\Gamma}, {\tilde{\pi}}^{\Gamma}_{*}p_{*}{{\underline{\mathbb{C}}}_{Y}}) & H^{n}(tot(\Lambda_{X_{\Gamma}}\mathscr{J}^{\bullet, *})) & ^{(c)}{{\tilde{E}}^{n,0}_{2}}=R^{n}(\Lambda_{X_{\Gamma}}{\tilde{\pi}}^{\Gamma}_{*})(p_{*}{{\underline{\mathbb{C}}}_{Y}})\\
     };
      \path (m-2-1) edge node[]{} (m-1-1);
      \path (m-3-1) edge node[]{} (m-2-1);
      \path (m-1-1) edge node[]{} (m-1-2);
      \path (m-2-1) edge node[]{} (m-2-2);
      \path (m-3-1) edge node[]{} (m-3-2);
      \path (m-3-2) edge node[]{} (m-2-2);
      \path (m-2-2) edge node[]{} (m-1-2);
      \path (m-1-3) edge node[below, label]{$\cong$} (m-1-2);
      \path (m-2-3) edge node[]{} (m-2-2);
      \path (m-3-3) edge node[above, label]{$\cong$} (m-3-2);
      \path (m-3-3) edge node[]{} (m-2-3);
      \path (m-2-3) edge node[]{} (m-1-3);
  \end{tikzpicture}
  \end{equation}}\\
Here the second row consists of the edge homomorphisms associated to the spectral sequence of the double complex $\Lambda_{M_{\Gamma}}{{\tilde{p}}^{-1}}\mathscr{J}^{\bullet,*}$ {\it(the first arrow is the edge hom. when filtered row wise and the second arrow is the edge hom. when filtered column wise and hence the notation (r), (c))}. Similarly the first and the third rows consist of the edge homomorphism associated to the GSS of the setups $C^{Y(\Gamma)}\overset{\pi^{\Gamma}_{*}}\rightarrow C^{M_{\Gamma}}\overset{\Lambda_{M_{\Gamma}}}\rightarrow Ab$ and $C^{X(\Gamma)}\overset{{\tilde{\pi}}^{\Gamma}_{*}}\rightarrow C^{X_{\Gamma}}\overset{\Lambda_{X_{\Gamma}}}\rightarrow Ab$ respectively. Indeed, the commutativity of each of the two squares enclosed by the first two rows and all the three columns  of (4) follows from the fact that we have a map of double complexes $\tilde{p}^{-1}\mathscr{J}^{\bullet, *}\rightarrow \mathscr{I}^{\bullet, *}$ extending the map of complexes ${\tilde{p}}^{-1}{\tilde{\pi}}^{\Gamma}_{*}\mathscr{J}^{\bullet}\overset{\delta^{\bullet}}\rightarrow \pi^{\Gamma}_{*}\mathscr{I}^{\bullet}$ in $C^{M_{\Gamma}}$, giving us a map between the associated row and column filtered complexes. The commutativity of each of the two squares enclosed by the second and the third rows and all the three columns of (4) follows from the fact that the natural transformation $\eta$ induces a map of double complexes $\Lambda_{X_{\Gamma}}\mathscr{J}^{\bullet, *}\rightarrow \Lambda_{M_{\Gamma}}{{\tilde{p}}^{-1}}\mathscr{J}^{\bullet,*}$ and hence between the associated row and column filtered complexes. The arrow obtained by composing the two arrows of the first column of (4) is the usual map induced by $\tilde{p}$ after making the obvious identification $\pi^{\Gamma}_{*}\underline{\mathbb{C}}_{Y}\cong \underline{\mathbb{C}}_{M_{\Gamma}}$. This is because the map ${\tilde{p}}^{-1}{\tilde{\pi}}^{\Gamma}_{*}(p_{*}{\underline{\mathbb{C}}}_{Y})\overset{\delta}\rightarrow \pi^{\Gamma}_{*}{\underline{\mathbb{C}}}_{Y}$ factors as ${\tilde{p}}^{-1}{\tilde{\pi}}^{\Gamma}_{*}(p_{*}{\underline{\mathbb{C}}}_{Y})\cong {\tilde{p}}^{-1}\tilde{p}_{*}\underline{\mathbb{C}}_{M_{\Gamma}}\rightarrow \underline{\mathbb{C}}_{M_{\Gamma}}\cong \pi^{\Gamma}_{*}\underline{\mathbb{C}}_{Y}$. The map $\mu$ is defined to be the composition of the two arrows in the third column of (4).
$\hfill$ $\Box$

Note that the map $\mu$ of lemma \ref{lemma2}, depends on the choice of the injective resolutions $0\rightarrow \underline{\mathbb{C}}_{Y}\rightarrow \mathscr{I}^{\bullet}$ in $C^{Y(\Gamma)}$, $0\rightarrow p_{*}\underline{\mathbb{C}}_{Y}\rightarrow \mathscr{J}^{\bullet}$ in $C^{X(\Gamma)}$ and the Cartan-Eilenberg resolutions $0\rightarrow \tilde{\pi}^{\Gamma}_{*}\mathscr{J}^{\bullet}\rightarrow \mathscr{J}^{\bullet,*}$ in $C^{X_{\Gamma}}$ and $0\rightarrow \pi^{\Gamma}_{*}\mathscr{I}^{\bullet}\rightarrow \mathscr{I}^{\bullet,*}$ in $C^{M_{\Gamma}}$. So we fix these choices for the proof of the lemma that follows.

GSS associated to the setup $C^{X(\Gamma)}\overset{\Lambda_{X}}\rightarrow \mathscr{M}\overset{F^{\Gamma}}\rightarrow Ab$, gives the edge homomorphism $H^{n}(\Gamma,p_{*}\underline{\mathbb{C}}_{Y}(X))\overset{v_{1}}\rightarrow R^{n}(\Lambda_{X_{\Gamma}}\tilde{\pi}^{\Gamma}_{*})(p_{*}\underline{\mathbb{C}}_{Y})$. Similarly GSS associated to the setup $C^{Y(\Gamma)}\overset{\Lambda_{Y}}\rightarrow \mathscr{M}\overset{F^{\Gamma}}\rightarrow Ab$, gives the edge homomorphism $H^{n}(\Gamma, \underline{\mathbb{C}}_{Y}(Y))\overset{v_{2}}\rightarrow R^{n}(\Lambda_{M_{\Gamma}}\pi^{\Gamma}_{*})(\underline{\mathbb{C}}_{Y})$.  
\begin{lemma}\label{lemma3}
Keeping the above notations, the following diagram commutes :
\begin{equation}
\begin{tikzpicture}[-stealth,
  label/.style = { font=\footnotesize }, baseline=(current  bounding  box.center)]
  \matrix (m)
    [
      matrix of math nodes,
      row sep    = 3.5em,
      column sep = 2.5em
    ]
    {
      H^{n}(\Gamma,p_{*}\underline{\mathbb{C}}_{Y}(X)) & H^{n}(\Gamma, \underline{\mathbb{C}}_{Y}(Y))\\      
      R^{n}(\Lambda_{X_{\Gamma}}\tilde{\pi}^{\Gamma}_{*})(p_{*}\underline{\mathbb{C}}_{Y}) & R^{n}(\Lambda_{M_{\Gamma}}\pi^{\Gamma}_{*})(\underline{\mathbb{C}}_{Y})\\
      };
    \path (m-1-1) edge node [above, label]{$id$} (m-1-2);
    \path (m-1-1) edge node [below, label]{$\cong$} (m-1-2);
    \path (m-1-1) edge node [right, label]{$v_{1}$} (m-2-1);
    \path (m-2-1) edge node [above, label]{$\mu$} (m-2-2);
    \path (m-1-2) edge node [right, label]{$v_{2}$} (m-2-2);
\end{tikzpicture}
\end{equation} 
\end{lemma}  
{\it Proof.}
We have the complexes $0\rightarrow\Lambda_{X}(p_{*}{\underline{\mathbb{C}}}_{Y})\rightarrow \Lambda_{X}{\mathscr{J}^{\bullet}}$ and $0\rightarrow \Lambda_{Y}{{\underline{\mathbb{C}}}_{Y}}\rightarrow \Lambda_{Y}\mathscr{I}^{\bullet}$ in $\mathscr{M}$ and choose their Cartan-Eilenberg resolutions $0\rightarrow\Lambda_{X}{\mathscr{J}^{\bullet}}\rightarrow \mathscr{K}^{\bullet, *}$ and $0\rightarrow\Lambda_{Y}\mathscr{I}^{\bullet}\rightarrow \mathscr{Z}^{\bullet, *}$ in $\mathscr{M}$ respectively. Since $p^{-1}$ is exact, as also mentioned in the proof of lemma \ref{lemma2}, there exists a morphism of complexes $p^{-1}\mathscr{J}^{\bullet}\overset{c^{\bullet}}\rightarrow \mathscr{I}^{\bullet}$ extending the obvious morphism $p^{-1}p_{*}{\underline{\mathbb{C}}}_{Y}\overset{c}\rightarrow \underline{\mathbb{C}}_{Y}$. There exists an obvious natural transformation $\beta : \Lambda_{X}\rightarrow \Lambda_{Y}p^{-1}$. We have the following commutative diagram :

\hspace{-2.5em}
\usetikzlibrary{matrix, calc}
\vspace{1.5em}
\hspace{10em}
\begin{tikzpicture}[-stealth,
  label/.style = { font=\footnotesize }]
  \matrix (m)
    [
      matrix of math nodes,
      row sep    = 2.5em,
      column sep = 3.5em
    ]
   {
\Lambda_{X}p_{*}{\underline{\mathbb{C}}}_{Y}=p_{*}{\underline{\mathbb{C}}}_{Y}(X) & \Lambda_{X}\mathscr{J}^{\bullet}\\
\Lambda_{Y}(p^{-1}p_{*}{\underline{\mathbb{C}}}_{Y})  & \Lambda_{Y}p^{-1}\mathscr{J}^{\bullet}\\
\Lambda_{Y}{\underline{\mathbb{C}}}_{Y}={\underline{\mathbb{C}}}_{Y}(Y) & \Lambda_{Y}\mathscr{I}^{\bullet}\\
     };
     \path (m-1-1) edge node[]{} (m-1-2);
     \path (m-2-1) edge node[]{} (m-2-2);
     \path (m-3-1) edge node[]{} (m-3-2);
     \path (m-1-1) edge node[right, label]{$\beta_{p_{*}{\underline{\mathbb{C}}_{Y}}}$} (m-2-1);
     \path (m-2-1) edge node[]{} (m-3-1);
     \path (m-1-2) edge node[right, label]{$\beta_{\mathscr{J}^{\bullet}}$} (m-2-2);
     \path (m-2-2) edge node[right, label]{$\Lambda_{Y}(c^{\bullet})=: \tilde{c}$} (m-3-2);
   \end{tikzpicture}\\
Thus we get a map of complexes $\Lambda_{X}\mathscr{J}^{\bullet}\overset{\Lambda_{Y}(c^{\bullet})\circ \beta_{\mathscr{J}^{\bullet}}}\longrightarrow \Lambda_{Y}\mathscr{I}^{\bullet}$ extending the identity morphism $\Lambda_{X}p_{*}{\underline{\mathbb{C}}}_{Y}\rightarrow \Lambda_{Y}(p^{-1}p_{*}{\underline{\mathbb{C}}}_{Y})\rightarrow \Lambda_{Y}{\underline{\mathbb{C}}}_{Y}$. Thus we get a map of Cartan-Eilenberg resolutions $\mathscr{K}^{\bullet, *}\rightarrow \mathscr{Z}^{\bullet, *}$ extending the map $\Lambda_{Y}(c^{\bullet})\circ \beta_{\mathscr{J}^{\bullet}}$ \cite [Ch. III, \S 7, Proposition 11 (b)]{gel-man}. 

On the other hand, $\eta$, $\zeta$ and $\beta$ induces the natural transformations $\eta^{'}: \Lambda_{X_{\Gamma}}{\tilde{\pi}}^{\Gamma}_{*}\rightarrow \Lambda_{M_{\Gamma}}{\tilde{p}}^{-1}{\tilde{\pi}}^{\Gamma}_{*}$, $\zeta^{'}: \Lambda_{M_{\Gamma}}{\tilde{p}}^{-1}{\tilde{\pi}}^{\Gamma}_{*}\rightarrow \Lambda_{M_{\Gamma}}\pi^{\Gamma}_{*}p^{-1}$ and $\beta^{'} : F^{\Gamma}\Lambda_{X}\rightarrow F^{\Gamma}\Lambda_{Y}p^{-1}$ respectively. Note that $\beta^{'}=\zeta^{'}\circ \eta^{'}$. 

Indeed we have the following commutative diagram, where the superscripts $(r)$, $(c)$ are as explained before :

\usetikzlibrary{matrix, calc}
{\tiny
\begin{equation}
\begin{tikzpicture}[-stealth,
  label/.style = { font=\footnotesize }, baseline=(current  bounding  box.center)]
  \matrix (m)
    [
      matrix of math nodes,
      row sep    = 5em,
      column sep = 2.1em
    ]
   {
H^{n}(\Gamma, p_{*}{{\underline{\mathbb{C}}}_{Y}}(X))= {^{(r)}{E^{''}}^{n,0}_{2}} & & & H^{n}(\Gamma, {\underline{\mathbb{C}}}_{Y}(Y))= {^{(r)}{E^{'}}^{n,0}_{2}}\\
H^{n}(tot(F^{\Gamma}\mathscr{K}^{\bullet,*})) & & & H^{n}(tot(F^{\Gamma}\mathscr{Z}^{\bullet,*}))\\
H^{n}(F^{\Gamma}\Lambda_{X}\mathscr{J}^{\bullet})= {^{(c)}{E^{''}}}^{n,0}_{2} & H^{n}(\Lambda_{M_{\Gamma}}{\tilde{p}}^{-1}{{\tilde{\pi}}^{\Gamma}_{*}}\mathscr{J}^{\bullet}) & H^{n}(F^{\Gamma}\Lambda_{Y}p^{-1}{\mathscr{J}^{\bullet}}) & H^{n}(F^{\Gamma}\Lambda_{Y}\mathscr{I}^{\bullet})={^{(c)}{E^{'}}^{n,0}_{2}}\\
      };
      \path (m-1-1) edge node[above, label]{$id$} (m-1-4);
      \path (m-2-1) edge node[]{} (m-2-4);
      \path (m-3-1) edge node[above, label]{via $\eta^{'}$} (m-3-2);
      \path (m-3-2) edge node[above, label]{via $\zeta^{'}$} (m-3-3);
      \path (m-3-3) edge node[above, label]{via $\tilde{c}$} (m-3-4);
      \path (m-1-1) edge node[below, label]{$\cong$} (m-1-4);
      \path (m-1-1) edge node[]{} (m-2-1);
      \path (m-3-1) edge node[right, label]{$\cong$} (m-2-1);
      \path (m-3-4) edge node[left, label]{$\cong$} (m-2-4);
      \path (m-1-4) edge node[]{} (m-2-4);
\end{tikzpicture}
\end{equation}}

In (6), the left and right columns are edge homomorphisms when GSS is applied to the set ups $C^{X(\Gamma)}\overset{\Lambda_{X}}\rightarrow \mathscr{M}\overset{F^{\Gamma}}\rightarrow Ab$, $C^{Y(\Gamma)}\overset{\Lambda_{Y}}\rightarrow \mathscr{M}\overset{F^{\Gamma}}\rightarrow Ab$ respectively. Since $\beta^{'}=\zeta^{'}\circ \eta^{'}$, the composition of the three arrows of the last row in (6) is same as the composition of the two arrows of the last column in (4) which is by definition the map $\mu$. The identity map of $\underline{\mathbb{C}}_{Y}(Y)$ in $\mathscr{M}$ induces a map from the injective resolution $Z^{0}(\mathscr{K}^{\bullet, *})$ to the injective resolution $Z^{0}(\mathscr{Z}^{\bullet, *})$ of $\underline{\mathbb{C}}_{Y}(Y)$ in $\mathscr{M}$ (both these resolutions are determined by the index *), thus the isomorphism which it induces at the level of cohomology can be viewed as the identity map of group cohomologies. 
$\hfill$ $\Box$

{\it Proof of Theorem \ref{theorem1} :} Now we further assume that $H^{q}(X,\mathbb{C})=0$ for all $q\gneq 0$. We have the following diagram :
\usetikzlibrary{matrix,calc}
{\tiny
\begin{equation}
\begin{tikzpicture}[-stealth,
  label/.style = { font=\footnotesize }, baseline=(current  bounding  box.center)]
  \matrix (m)
    [
      matrix of math nodes,
      row sep    = 4.9em,
      column sep = 5em
    ]
   {
H^{n}(\Gamma, \underline{\mathbb{C}}_{X}(X)) & H^{n}(\Gamma, p_{*}\underline{\mathbb{C}}_{Y}(X)) & H^{n}(\Gamma, \underline{\mathbb{C}}_{Y}(Y)) & H^{n}(\Gamma, \mathcal{O}_{Y}(Y))\\  
R^{n}(\Lambda_{X_{\Gamma}}\tilde{\pi}^{\Gamma}_{*})(\underline{\mathbb{C}}_{X}) & R^{n}(\Lambda_{X_{\Gamma}}\tilde{\pi}^{\Gamma}_{*})(p_{*}\underline{\mathbb{C}}_{Y}) & R^{n}(\Lambda_{M_{\Gamma}}\pi^{\Gamma}_{*})(\underline{\mathbb{C}}_{Y}) & R^{n}(\Lambda_{M_{\Gamma}}\pi^{\Gamma}_{*})(\mathcal{O}_{Y})\\
H^{n}(X_{\Gamma}, \tilde{\pi}^{\Gamma}_{*}\underline{\mathbb{C}}_{X}) & H^{n}(X_{\Gamma}, \tilde{p}_{*}\underline{\mathbb{C}}_{M_{\Gamma}}) & H^{n}(M_{\Gamma},\underline{\mathbb{C}}_{M_{\Gamma}}) & H^{n}(M_{\Gamma},\mathcal{O}_{M_{\Gamma}})\\   
      };
      \path (m-1-1) edge node[above, label]{$\Lambda_{X}(t)_{*}$} (m-1-2);
      \path (m-1-1) edge node[below, label]{$\cong$} (m-1-2);
      \path (m-1-2) edge node[above, label]{$id$} (m-1-3);
      \path (m-1-2) edge node[below, label]{$\cong$} (m-1-3);
      \path (m-1-3) edge node[above, label]{$\Lambda_{Y}(\iota)_{*}$} (m-1-4);
      \path (m-2-1) edge node[above, label]{$t_{*}$} (m-2-2);
      \path (m-2-2) edge node[above, label]{$\mu$} (m-2-3);
      \path (m-2-3) edge node[above, label]{$\iota_{*}$} (m-2-4);
      \path (m-3-1) edge node[above, label]{$t_{*}$} (m-3-2);
      \path (m-3-2) edge node[above, label]{via $\tilde{p}$} (m-3-3);
      \path (m-3-3) edge node[above, label]{$\iota_{*}$} (m-3-4);
      \path (m-1-1) edge node[right, label]{$\cong$} (m-2-1);
      \path (m-1-1) edge node[left, label]{$u_{1}$} (m-2-1);
      \path (m-1-2) edge node[right, label]{$\cong$} (m-2-2);
      \path (m-1-2) edge node[left, label]{$v_{1}$} (m-2-2);      
      \path (m-1-3) edge node[left, label]{$v_{2}$} (m-2-3);
      \path (m-1-4) edge node[right, label]{$v$} (m-2-4);
      \path (m-3-1) edge node[below, label]{$\cong$} (m-3-2);
      \path (m-2-1) edge node[below, label]{$\cong$} (m-2-2);      
      \path (m-3-1) edge node[left, label]{$$} (m-2-1);
      \path (m-3-1) edge node[left, label]{$s$} (m-2-1);
      \path (m-3-2) edge node[left, label]{$e_{1}$} (m-2-2);
      \path (m-3-3) edge node[right, label]{$\cong$} (m-2-3);
      \path (m-3-3) edge node[left, label]{$e_{2}$} (m-2-3);
      \path (m-3-4) edge node[right, label]{$\cong$} (m-2-4);
      \path (m-3-4) edge node[left, label]{$w$} (m-2-4);
      \end{tikzpicture}
\end{equation}}\\
The third vertical arrow arising from the last row is an isomorphism because $\pi$ is a (holomorphic) covering map \cite [Ch. 5, Corollary 1 after Theorem 5.3.1]{grothendieck}. The top square enclosed by all the three rows and the second and the third columns of (7) is same as diagram (5) which commutes by lemma \ref{lemma3}, and the bottom square enclosed by all the three rows and the second and the third columns of (7) is same as diagram (3) which commutes by lemma \ref{lemma2}. A simple argument using the functorial properties of spectral sequences of filtered complexes accounts for the commutativity of the rest of the squares of diagram (7). Thus (7) commutes and this proves the commutativity of (2). $\hfill$ $\Box$

Now let $G$ be a connected complex semisimple Lie group, $\Gamma$ be a cocompact, irreducible and torsionless lattice in $G$ and $K$ be a maximal compact subgroup of $G$. Assume $\Gamma$ acts by left multiplication and $K$ acts by right multiplication on $G$. Taking $Y=G$, $\Gamma$= the lattice $\Gamma$, $K=$ the maximal compact $K$ and considering the obvious left action of $\Gamma$ on $Y$ (which is free and proper, and hence a covering space action), the obvious right action of $K$ on $Y$ (which is smooth free and proper), and $X=G/K$ (contractible smooth manifold), we see that indeed the hypotheses of theorem \ref{theorem1} are satisfied. Moreover in this case, the usual left action of $\Gamma$ on $X$ is free and proper, forcing $X_{\Gamma}$ to be a smooth manifold and $\tilde{\pi}$ to be a smooth covering. We note that $\Gamma\backslash Y = M_{\Gamma}$ is a compact complex manifold and $X_{\Gamma}\cong\Gamma\backslash G/K\cong M_{\Gamma}/K$ is also a smooth compact manifold. 

{\it Proof of Theorem \ref{t1} : }The map $\psi$ obtained by applying theorem \ref{theorem1} to $Y$, $\Gamma$ and $K$ as in the preceding paragraph, is an isomorphism since $\tilde{\pi}$ being a smooth covering, $s$ is an isomorphism \cite [Ch. 5, Corollary 1 after Theorem 5.3.1]{grothendieck}, $G$ being Stein, $H^{q}(G,\mathcal{O}_{G})=0$ for all $q>0$ \cite[Ch. 4, \S 2, Theorem B]{onishchik} forcing $v$ to be an isomorphism \cite [Ch. 5, Proposition 5.2.5]{grothendieck} and $\Lambda_{Y}(\iota)_{*}$ is an isomorphism \cite[proof of Theorem 1]{akhiezer}.   $\hfill$ $\Box$
\begin{remark}
It is easy to see that the isomorphism in \cite[Theorem 1]{akhiezer} when $p=0$, is same as $\psi$ as obtained in the proof of theorem \ref{t1}.
\end{remark}
\section{Applications to Picard Groups}
Let $G$ be a simply connected complex semisimple Lie group, $\Gamma$ be a cocompact, irreducible and torsionless lattice in $G$ and $K$ be a maximal compact subgroup of $G$. Assume $\Gamma$ acts by left multiplication and $K$ acts by right multiplication on $G$. Let $M_{\Gamma}= \Gamma\backslash G$, $X=G/K$ and $X_{\Gamma}=\Gamma\backslash X$.

In this section we compute the Picard group of $\Gamma\backslash G$ with rank($G$) $=1$ and $2$. In other words we provide a proof of theorem \ref{t2}.

The classification theorem says that the only simply connected complex semisimple Lie group of rank one is $SL(2,\mathbb{C})$ and those of rank 2 are $SL(3,\mathbb{C})$, $SL(2,\mathbb{C})\times SL(2, \mathbb{C})$, $Spin(5)$ and the exceptional group $G_{2}$. 

Applying Leray-Serre spectral sequence to the smooth principal $K-$bundle $\tilde{p}:M_{\Gamma}\rightarrow X_{\Gamma}$ and using $H^{q}(K,\mathbb{Z})=0$ for $q=1,2$ (since $K$ is compact, simply connected and semisimple), we get that $\tilde{p}^{*}:H^{q}(X_{\Gamma},\mathbb{Z})\rightarrow H^{q}(M_{\Gamma},\mathbb{Z})$ is an isomorphism for $q=1,2$. Same conclusion holds with $\mathbb{C}$-coefficients.
 
{\it Proof of Theorem \ref{t2} : }Let $A$ be the torsion subgroup of $H^{2}(M_{\Gamma},\mathbb{Z})$. If rank($G$) $=2$, then $H^{1}(X_{\Gamma},\mathbb{C})=0$ \cite[Corollary 6]{picardgroups}. Thus $H^{1}(M_{\Gamma},\mathbb{C})=0$ forcing $H^{1}(M_{\Gamma},\mathbb{Z})=0$ (since $H_{1}(M_{\Gamma},\mathbb{Z})$ is finite) and $H^{1}(M_{\Gamma},\mathcal{O}_{M_{\Gamma}})=0$ \cite[Theorem 1]{akhiezer}. Long exact sequence in sheaf cohomology induced by the exponential exact sequence of sheaves of abelian groups on $M_{\Gamma}$ : 
\[0\rightarrow \underline{\mathbb{Z}}_{M_{\Gamma}}\rightarrow \mathcal{O}_{M_{\Gamma}}\rightarrow \mathcal{O}_{M_{\Gamma}}^{*}\rightarrow 0\]

\noindent implies that the first Chern class map $c_{1}: Pic(M_{\Gamma})\rightarrow H^{2}(M_{\Gamma},\mathbb{Z})$ is injective. Since image of $c_{1}$ is same as kernel of the map $i:H^{2}(M_{\Gamma},\mathbb{Z})\rightarrow H^{2}(M_{\Gamma},\mathcal{O}_{M_{\Gamma}})$ induced by the inclusion of the constant sheaf $\mathbb{Z}$ on $M_{\Gamma}$ into $\mathcal{O}_{M_{\Gamma}}$, we only need to show that the above kernel is same as the kernel of the map $i_{1}:H^{2}(M_{\Gamma},\mathbb{Z})\rightarrow H^{2}(M_{\Gamma},\mathbb{C})$ induced by the inclusion of the constant sheaf $\mathbb{Z}$ into the constant sheaf $\mathbb{C}$ on $M_{\Gamma}$ (which is same as $A$). Let $i_{2}:H^{2}(M_{\Gamma},\mathbb{C})\rightarrow H^{2}(M_{\Gamma},\mathcal{O}_{M_{\Gamma}})$ be the map induced by the inclusion of the constant sheaf $\mathbb{C}$ on $M_{\Gamma}$ into $\mathcal{O}_{M_{\Gamma}}$. Note that $i$ factors as $i=i_{2}\circ i_{1}$. Since $\tilde{p}^{*}:H^{2}(X_{\Gamma},\mathbb{C})\rightarrow H^{2}(M_{\Gamma},\mathbb{C})$ is an isomorphism, theorem \ref{t1} implies that $i_{2}$ is an isomorphism. Thus $ker$ $i=$ $ker$ $i_{1}$. This proves the theorem in the rank two case.

Now let rank($G$) $=1$. Same arguments as in the rank two case show that the image of the map $c_{1}:Pic(M_{\Gamma})\rightarrow H^{2}(M_{\Gamma},\mathbb{Z})$ is $A$, but unlike rank two case $Pic^{0}(M_{\Gamma})$ need not vanish. In fact $Pic^{0}(M_{\Gamma})=\mathbb{C}^{r}/\mathbb{Z}^{r}$, where $r=$ rank($\Gamma/[\Gamma,\Gamma]$) \cite [Theorem 1]{picardgroups}. The rank one case now follows from the fact $Pic(M_{\Gamma})=Pic^{0}(M_{\Gamma})\oplus im(c_{1})$. $\hfill$ $\Box$\\

{\bf Acknowledgements :} I am indebted to P. Sankaran for his striking claims which led to this research work. I would like to thank him for his suggestions which helped me to improve the exposition. I also take this oppurtunity to thank V. Balaji and D.S. Nagaraj for their valuable comments.

\end{document}